\documentclass[12pt]{article}
\usepackage[utf8]{inputenc}
\usepackage{amsmath}
\usepackage{amsfonts}
\usepackage{amssymb}
\usepackage{amsthm}
\usepackage{relsize}
\usepackage{array}
\usepackage{multirow}
\usepackage{tabu}
\usepackage{soul}
\usepackage{graphicx}
\usepackage{epigraph}
\usepackage{float}
\usepackage[english]{babel}
\usepackage[usenames, dvipsnames]{color}
\usepackage{fancyhdr}
\usepackage[Sonny]{fncychap}
\usepackage{tikz-cd}
\usepackage{tkz-euclide}
\usepackage{sseq}
\usepackage{newfloat}
\usepackage{amsmath,mathtools}
\usepackage{pgfplots}
\usepackage{geometry}
 \usepackage{sectsty}
\sectionfont{\fontsize{14}{15}\selectfont}

\usepackage[all]{xy}
\DeclareFloatingEnvironment[fileext=dia]{diagram}

\theoremstyle{definition}

\newtheorem{remark}{Remark}

\usepackage[most]{tcolorbox}

\makeatletter
\newtcolorbox{note}[1][]{%
	breakable,
	enhanced jigsaw, 
	borderline west={3pt}{0pt}{black!10!white}, 
	borderline south={1pt}{0pt}{black!10!white}, 
	borderline east={1pt}{0pt}{black!10!white},
	borderline north={1pt}{0pt}{black!10!white},
	sharp corners, 
	boxrule=0pt, 
	attach title to upper, 
	left=0pt,
	right=0pt,
	top=0pt,
	bottom=0pt,
	boxsep=5pt,
	colback=white,
	frame hidden,
	#1
}
\makeatother

\makeatletter
\newtcolorbox{note1}[1][]{%
	breakable,
	enhanced jigsaw, 
	sharp corners, 
	boxrule=0pt, 
	attach title to upper, 
	fontupper=\linespread{1.1}\fontfamily{qpl}\selectfont,
	fontlower=\linespread{1.1}\fontfamily{qpl}\selectfont, 
	left=0pt,
	right=0pt,
	top=0pt,
	bottom=0pt,
	boxsep=3pt,
	colback=green!3!white,
	frame hidden,
	before skip=10pt plus 2pt,after skip=10pt plus 2pt,
	#1
}
\makeatother

\makeatletter
\newcommand\tabfill[1]{%
	\dimen@\linewidth
	\advance\dimen@\@totalleftmargin
	\advance\dimen@-\dimen\@curtab
	\parbox[t]\dimen@{#1\ifhmode\strut\fi}%
}
\makeatother

\usepackage{tabularx}
\usepackage{imakeidx}
\usepackage{hyperref}
\hypersetup{colorlinks={true},linkcolor={blue},citecolor={blue},urlcolor={blue}, linktoc=all}  
 
 \usepackage[capitalize, nameinlink]{cleveref}
 \crefdefaultlabelformat{#2\textbf{#1}#3}
 \crefname{figure}{Figure}{Figures} 
 \Crefname{figure}{Figure}{Figures}
 \crefname{table}{Table}{Tables}
 \Crefname{table}{Table}{Tables}
 \crefname{section}{\S\hspace{-1mm}}{\S\hspace{-1mm}}
 \Crefname{section}{\S\hspace{-1mm}}{\S\hspace{-1mm}}
 \crefname{equation}{}{}
 \Crefname{equation}{}{}
 \crefname{example}{Geometric Pattern}{Geometric Patterns} 
 \Crefname{example}{Geometric Pattern}{Geometric Patterns}

 \usepackage{graphicx,tipa}

 \usepackage{footnote}

\newcommand{\R}{\mathbb{R}}

\begin{document}

\title{\textbf{Systems of Equations in Elamite Mathematics}}

\author{Nasser Heydari\footnote{Email: nasser.heydari@mun.ca}~ and  Kazuo Muroi\footnote{Email: edubakazuo@ac.auone-net.jp}}

\maketitle

\begin{abstract}
This article studies the systems  of equations appearing  in the Susa Mathematical Texts (\textbf{SMT}) and the  different approaches used by the Susa scribes to solve them.      
\end{abstract}

\section{Introduction}
Similar to   quadratic equations, systems of equations  appear in several texts of the Susa mathematical texts including  \textbf{SMT No.\,8}, \textbf{SMT No.\,11}, \textbf{SMT No.\,17},  \textbf{SMT No.\,18}, and \textbf{SMT No.\,19}. These texts belong to   26 clay tablets excavated from Susa in  southwest Iran by French archaeologists in 1933. The texts of all the Susa mathematical texts (\textbf{SMT}) along with their interpretations were first published in 1961 under the title \textit{Textes Math\a'{e}matiques de Suse}, which is known as the \textbf{TMS} in the literature (see \cite{BR61}).

The upper and lower parts of \textbf{SMT No.\,8}\footnote{The reader can see   this tablet on the website of the Louvre's collection. Please see \url{https://collections.louvre.fr/en/ark:/53355/cl010186534} for obverse  and   reverse, as well as \url{https://collections.louvre.fr/en/ark:/53355/cl010186431} for its side.}      are lost, and  the reverse of the tablet is   almost completely destroyed. Fortunately, however, two algebraic problems regarding simultaneous equations are preserved on the obverse. In each problem, the   corresponding  system  of equations  is  changed   to  a quadratic equation  by using a new variable and finally the obtained quadratic equation  is solved by  completing the square.   As far as we can judge, two more problems should have been on the obverse. 

 \textbf{SMT No.\,11}\footnote{The reader can see   this tablet on the website of the Louvre's collection. Please see \url{https://collections.louvre.fr/en/ark:/53355/cl010186537} for obverse  and   reverse.}    contains two problems one of which can be described as a primitive  indeterminate equation  and the other as a complicated system of quadratic equations. Although   some parts of both   problems are unintelligible to us, especially in the statement of the second problem, there is no doubt about their mathematical meanings.  At the end of the first problem in line 7,   the scribe of the tablet   stops writing and leaves the solution incomplete. This leaves a very  noticeable blank space on the obverse of the tablet.

\textbf{SMT No.\,17}\footnote{The reader can see   this tablet on the website of the Louvre's collection. Please see \url{https://collections.louvre.fr/en/ark:/53355/cl010186427} for obverse  and   reverse.}    contains, pace the \textbf{TMS}, only one problem which is reduced to a typical example of the simultaneous quadratic equations found in  Babylonian mathematics. In the solution the scribe does not use his routine method, that is, making use of symmetric expressions $x+y = a$  and $xy = b$, but rather uses the factorization method called {\fontfamily{qpl}\selectfont \textit{mak\c{s}arum}}. In the latter half of the reverse of the tablet he confirms the answers obtained to be correct.

 \begin{remark}
 	Although a system of equations appears in \textbf{SMT No.\,18}, we do not consider it  here and will discuss it in another  paper on the similarity of triangles and intercept theorem.  
 \end{remark}
 
The text  of \textbf{SMT No.\,19}\footnote{The reader can see   this tablet on the website of the Louvre's collection. Please see \url{https://collections.louvre.fr/en/ark:/53355/cl010186429} for obverse  and   reverse.}     contains two problems, one on the obverse and the other on the reverse of the tablet, both of which deal with simultaneous equations  concerning Pythagorean triples (see \cite{HM23-3}, for a discussion on Pythagorean theorem and triples in the \textbf{SMT}). In the first problem the diagonal  (of a rectangle) or the hypotenuse (of a right triangle) is called {\fontfamily{qpl}\selectfont tab}  ``friend, partner'' which reminds us of the fact that a Pythagorean triple was called {\fontfamily{qpl}\selectfont illat} ``group, clan'' in Babylonian mathematics. The scribe of this tablet handles these equations skillfully, especially in the second problem which might  be one of the most complicated systems of simultaneous equations in  Babylonian mathematics.

\section{Systems of Equations}
A system of equations is a finite set of simultaneous equations with finite unknown variables. To be more precise, let $x_1,x_2,\cdots,x_n$ be $n$ variables and $f_1(x_1,x_2,\cdots,x_n)$, $f_2(x_1,x_2,\cdots,x_n)$, $\ldots$, $f_m(x_1,x_2,\cdots,x_n)$ be $m$ functions of $x_i$. The set of equations 
\begin{equation}\label{eq-aa}
	\begin{cases}
		f_1(x_1,x_2,\ldots,x_n)=0 \\
		f_2(x_1,x_2,\ldots,x_n)=0 \\
		~~~~~~~~~~~\vdots \\
		f_m(x_1,x_2,\ldots,x_n)=0 \\
	\end{cases}
\end{equation} 
is a \textit{system of simultaneous equations} with $m$ equations and $n$  unknowns. The functions $f_j$ can be any  algebraic functions with respect to  variables $x_i$. For example, 

\begin{equation}\label{eq-ab}
	\begin{cases}
		x^2-yx^3z-\sin(xz)=0 \\
		\log(x^3-y)+\sqrt{z+4}=0 \\
		xyz-z^4+y^3-5=0  
	\end{cases}
\end{equation} 
is a system of three equations   and three unknowns $x,y,z$. 

Any solution of \cref{eq-ab} simultaneously satisfies  $m$ polynomial equations. So, the solution set of  \cref{eq-ab} is the intersection of all zero sets of   polynomials $f_1,f_2,\ldots,f_n$. It should be noted that finding the    solutions of  systems of equations is not possible in general. However by using the many numerical algorithms at our disposal we can find good estimates of solutions (see \cite{BHSW}). 
  
\section{Systems of Polynomial Equations}
A particular class of systems of equations is the one whose functions are all polynomials. Such   systems are    usually called  \textit{systems of polynomial equations}. 

Recall that a polynomial function of  $n$  real variables $x_1,x_2,\ldots,x_n$ is an algebraic expression  like
\[ p(x_1,x_2,\ldots,x_n)=\sum_{0\leq i_1,i_2,\ldots,i_n\leq n}a_{i_1,i_2,\ldots,i_n}   x^{i_1}_{1}x^{i_2}_{2}\ldots x^{i_n}_{n} \]
 where $a_{i_1,i_2,\ldots,i_n}$ are real numbers.  The maximum value of summations $i_1+i_2+\cdots+i_n$ for which $a_{i_1,i_2,\ldots,i_n}  \neq 0$ is the \textit{degree} of  $p$ and denoted by $\deg(p)$.    The equation $p(x_1,x_2,\ldots,x_n)=0$ is called a \textit{polynomial equation}. The solution of such an equation is the set of all $n$-tuples $(z_1,z_2,\ldots,z_n)\in \R^n$ such that 
 \[  \sum_{0\leq i_1,i_2,\ldots,i_n\leq n}a_{i_1,i_2,\ldots,i_n}   z^{i_1}_{1}z^{i_2}_{2}\ldots z^{i_n}_{n}=0. \]
 This set is usually called  the \textit{zero set} of $p$, denoted by $Z_p$, which is a subset of the $n$-dimensional Euclidean spare $\R^n$. Such zero sets are also   called   \textit{real algebraic curves} in algebraic geometry.  For example, the zero set of cubic polynomial function $p(x,y)=(x-1)(x^2+y^2)+4x^2$ is a curve in the plane known as the \textit{Conchoid of de Sluze}.

 \begin{figure}[H]
 	\centering
 	\includegraphics[scale=1]{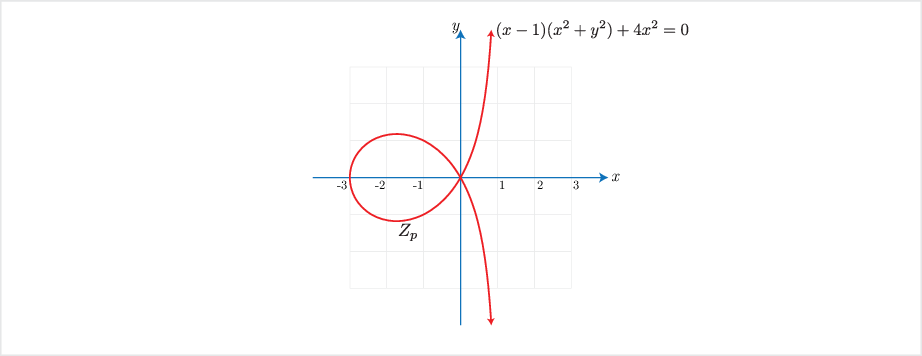}
 	\caption{A real algebraic curve: the Conchoid of de Sluze}
 	\label{Figure1}
 \end{figure}

\section{Methods to Solve a System of Equations}
Solving a system of equations is not easy in general and one needs to use different approaches and techniques to find the  possible solutions. The choice of appropriate techniques to solve a system of equations  depends on the formulas of its equations and the number of its unknowns. A system of equations  might have no solution, have a finite number of solutions or an infinite number of solutions.  
 
 \subsection{Gaussian Elimination: Linear Equations}
If the polynomial equations in a system of equations are all \textit{linear} polynomials, i.e., their degree is 1, then the system is called a \textit{system of linear equations}. The general form of a system of linear equations is as follows:
\begin{equation}\label{eq-bb}
	\begin{cases}
		a_{11}x_1+\ldots+a_{1n}x_n =b_1 \\
		a_{21}x_1+\ldots+a_{2n}x_n =b_2 \\ 
		~~~~~~~~~~~~\vdots\\
		a_{m1}x_1+\ldots+a_{mn}x_n =b_m \\
	\end{cases}
\end{equation} 
where $a_{ij}$ and $b_{i}$ are real numbers and $x_{j}$ are unknown variables.

For example, consider  
\begin{equation}\label{eq-ba}
	\begin{cases}
		x+y-z=1 \\
		x-y+z=1  
	\end{cases}
\end{equation} 
which is a system of two linear  equations and three variables. It can be shown  that the solution of this system is the set of triples 
\[ L = \{(1,t,t) ~|~ t\in \R \}. \]

\begin{figure}[H]
	\centering
	\includegraphics[scale=1]{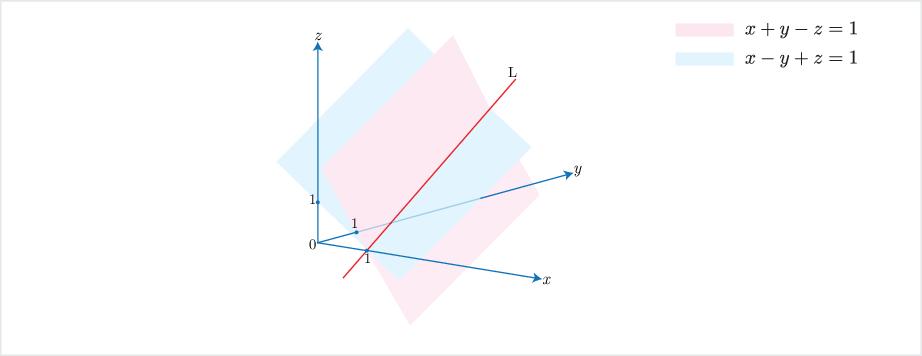}
	\caption{Solution of a system of linear equations}
	\label{Figure2}
\end{figure}

This set  is a line in the three-dimensional space  passing through $(1,0,0)$ and  generated by the vector $(0,1,1)$. \cref{Figure2} shows the solution of the system as the intersection of two planes in the space, which is the red-colored  line $L$.

Systems of linear equations  are usually solved by  a technique  from linear algebra known  as the \textit{row reduction} or \textit{Gaussian elimination} (see \cite{HK71}, for more details).  In fact, if we form the    matrices 
\[A= \begin{bmatrix}
	a_{11} &a_{12} & \cdots & a_{1n}\\
	a_{21} &a_{22} & \cdots & a_{2n}\\
	\vdots  & \vdots & \ddots&   \vdots \\
	a_{m1} &a_{m2} & \cdots & a_{mn}\\
\end{bmatrix}, ~X=\begin{bmatrix}
x_{1} \\
x_{2}\\
\vdots  \\
x_{n} \\
\end{bmatrix}, ~B= \begin{bmatrix}
	b_{1} \\
	b_{2}\\
	\vdots  \\
	b_{m} \\
\end{bmatrix} \]
then we can rewrite the system of \cref{eq-bb} as the matrix equation $AX=B$. In this case, $A$ is usually  called the \textit{coefficient matrix} of the system.  

The essence of the  method to solve  the matrix equation $AX=B$ is to apply three \textit{elementary operations} to the rows  of the  coefficient matrix $A$ and matrix $B$ simultaneously in order to reduce $A$  to an upper triangular matrix.   For example, the system of linear equations
\begin{equation}\label{eq-b}
	\begin{cases}
		2x-y+3z=2 \\
		~~x+y+2z=1 \\
	-x+2y+z=0 \\
	\end{cases}
\end{equation}  
can be reduced to the following reduced system by doing some reduction operations on the rows of its coefficient matrix:
\begin{equation}\label{eq-c}
	\begin{cases}
		x+y+4z=2 \\
		~~~~~~~y+z=\frac{1}{3} \\
		~~~~~~~~~~~~z=\frac{1}{2} \\
	\end{cases}
\end{equation}  
This reduced system   easily gives the solution  as $x=\frac{1}{6}, y=-\frac{1}{6},$ and  $z=\frac{1}{2}$. 

 \subsection{Elimination of Variables}
The simplest method to solve a system of equations might be the \textit{elimination of variables} in which one repeatedly eliminates variables to get a system  with less variables. If these equations provide  us with a solution, we can substitute it in the main system to get a new system with less variables. By repeating the same pattern, ultimately, one may be able to find the solution of the main system. For example, consider the system of   equations
\begin{equation}\label{eq-d}
	\begin{cases}
		z+xy+z^2=3 \\
		y^2-xz+z^2=1 \\
		xy+z=2 \\
	\end{cases}
\end{equation} 
If we substitute  $xy=2-z$ in the first equation and simplify, we get  the quadratic equation $ z^2 =1$ which gives $z=1$ and $z=-1$. Using $z=1$ in the third equation, we obtain $xy=1$ and so $x=\frac{1}{y}$. If we set both  $z=1$ and $x=\frac{1}{y}$    in the second equation, we get $y^2 = \frac{1}{y}$ which gives $y^3=1$ or $y=1$. So, the first solution is the triple $(1,1,1)$. In a similar fashion, for the case $z=-1$, we get $x=-\sqrt[3]{9}$ and $y=-\sqrt[3]{3}$. Thus, the second solution is the triple $(-\sqrt[3]{9},-\sqrt[3]{3},-1)$.

 \subsection{Change of Variables}
Another method to solve a system of equations is the change of variable. Generally speaking,  a change of variables is a basic technique  used to simplify problems in which the original variables are replaced with functions of those variables. This might cause  the problem to  become simpler, or equivalent to a better understood problem. This method has been applied in algebra to solve algebraic equations for a   long time and can be considered as one of the first techniques  invented by mathematicians   for solving equations. Although the first appearance of this method seems to be in algebraic equations, there are   applications of it in other branches of mathematics such as integration, differentiation, and differential equations.   

One category of systems of equations with two variables $x,y$  appearing in Babylonian texts contains those systems whose equations involve the sum and product of the two variables. In such systems, the new variables $u=x+y$ and $v=xy$ are appropriate choices for new variables and  one might get a simpler system of equations by substituting  the values of $u$ and $v$. For each pair of values of $u,v$, one needs to solve  the  quadratic equation $z^2-uz+v=0$. For instance, the system 
\begin{equation}\label{eq-e}
	\begin{cases}
		x^2-xy+y^2=4 \\
		x^2y+xy^2=16 \\
	\end{cases}
\end{equation} 
can be solved by using  the new variables $u=x+y$ and $v=xy$. Since $x^2-xy+y^2= (x+y)^2-3xy$ and $x^2y+xy^2=xy(x+y)$, we get a new system
 \begin{equation}\label{eq-e}
 	\begin{cases}
 		u^2-3v=4 \\
 		uv=16 \\
 	\end{cases}
 \end{equation} 
 Setting $v=\frac{16}{u}$ in the first equation gives $u^3-4u-48=0$. Since $u^3+4u-48=(u-4)(u^2+4u+12)$, we get only $u=4$ implying that $v=4$. So, we need to solve the quadratic equation $z^2-4z+4=0$ whose solutions can be found by the quadratic formula or completing the square. Ultimately, we get $x=y=2$.

\section{Systems of Equations in the  SMT} 
 In this part, we investigate the systems of equations found in  some texts  of the   \textbf{SMT} and set out our mathematical interpretations. 

 In order to understand these texts   correctly, it is helpful to introduce a practice of Babylonian mathematics in solving equations. Since the Babylonians did not have any mathematical symbols, they sometimes used, in the intermediate calculations, the modifiers {\fontfamily{qpl}\selectfont \textit{sarrum} (= lul)} ``provisional, false'' and {\fontfamily{qpl}\selectfont \textit{k\={\i}num} (= gi-na)} ``true'' concerning {\fontfamily{qpl}\selectfont u\v{s}} ``length'', {\fontfamily{qpl}\selectfont sag} ``width'', {\fontfamily{qpl}\selectfont a-\v{s}\`{a}} ``area'', etc. Although these modifiers are not used in these texts by the Susa scribes, we have put   them in round brackets in the following translation for   reader's convenience.  

\subsection{SMT No.\,8}

\subsubsection*{First Problem}\label{SS-P1-SMT8} 
By checking the text carefully, one can  notice that a part of the double-ruling line is left above line 1  which implies that our ``first'' problem was preceded by another problem.

\subsubsection*{Transliteration of First Problem}\label{SSS-P1TI-SMT8} 
\begin{note1} 
	\underline{Obverse:  Lines 1-10}  \\
	(L1)\hspace{1mm} [a-\v{s}\`{a} 10 4-\textit{at} sag \textit{a-na} sag dah] \textit{a-na} 3 \textit{a-li}-[\textit{ik a-na-di} sag ugu u\v{s}]\\
	(L2)\hspace{1mm} [5 dir]ig za-e [4 \textit{r}]\textit{e-ba-ti ki-ma} sag gar \textit{re-ba}-[\textit{at} 4 \textit{le-q\'{e}} 1 \textit{ta-mar}]\\
	(L3)\hspace{1mm} [1 \textit{a-na}] 3 \textit{a-li-ik} 3 \textit{ta-mar} 4 \textit{re-ba-at} sag \textit{a-na} 3 d[ah 7 \textit{ta-mar}]\\
	(L4)\hspace{1mm} 7 \textit{ki-ma} u\v{s} gar 5 dirig \textit{a-na} \textit{na-s\'{i}-ih} u\v{s} gar 7 u\v{s} \textit{a-na} 4 [sag \textit{i-\v{s}\'{i}}]\\
	(L5)\hspace{1mm} 28 \textit{ta-mar} 28 a-\v{s}\`{a} 28 \textit{a-na} 10 a-\v{s}\`{a} \textit{i-\v{s}\'{i}} 4,40 \textit{ta-mar}\\
	(L6)\hspace{1mm} [5] \textit{na-s\'{i}-ih} u\v{s} \textit{a-na} 4 sag \textit{i-\v{s}\'{i}} 20 \textit{ta-mar} 1/2 \textit{he-pe} 10 \textit{ta-mar} 10 nigin\\
	(L7)\hspace{1mm} [1,40] \textit{ta-mar} 1,40 \textit{a-na} 4,40 dah 4,41,40 \textit{ta-mar mi-na} \'{i}b-si 2,10 \textit{ta-mar}\\
	(L8)\hspace{1mm} [10 \textit{le}]-\textit{q\'{i}} \textit{a-na} 2,10 dah 2,20 \textit{ta-mar mi-na a-na} 28 a-\v{s}\`{a} gar \textit{\v{s}\`{a}} 2,20 \textit{i-na}-[\textit{di-n}]\textit{a}\\
	(L9)\hspace{1mm} [5 gar] 5 \textit{a-na} 7 \textit{i-\v{s}\'{i}} 35 \textit{ta-mar} 5 \textit{na-s\'{i}-ih} u\v{s} \textit{i-na} 35 zi\\
	(L10)\hspace{-1mm} [30 \textit{ta}]-\textit{mar} 30 u\v{s} 5 u\v{s}(sic) \textit{a-na} 4 sag \textit{i-\v{s}\'{i}} 20 \textit{ta-mar} 20 u\v{s}(sic) 
\end{note1}

\subsubsection*{Translation of First Problem}\label{SS-P1TR-SMT8} 
\underline{Obverse:  Lines 1-10} 
\begin{tabbing}
	\hspace{12mm} \= \kill
	(L1)\> \tabfill{The area is 10,0. Add one fourth of the width to the width. Multiply (one fourth   of the width) by 3. I lay down that the (provisional) widths exceed the length  by 5.}\\
	(L2)\>  \tabfill{You, put down 4 of one fourth as the width. Take one fourth of 4, (and) you see   1.}\\ 
	(L3)\> \tabfill{Multiply 1 by 3, (and) you see 3. Add 4 of one fourth of the width to 3, (and)  you see 7.} \\
	(L4)\> \tabfill{Put down 7 as the (provisional) length. Put down the excess 5 for that which   is subtracted from the (provisional) length. Multiply the (provisional) length 7   by 4 of the width, (and)}\\
	(L5)\> \tabfill{you see 28. 28 is the (coefficient of provisional) area. Multiply 28 by the area   10,0 (and) you see 4,40,0.} \\
	(L6)\> \tabfill{Multiply 5, which is subtracted from the (provisional) length, by 4 of the width,  (and) you see 20. Halve (it, and) you see 10. Square 10, (and)}\\
	(L7)\> \tabfill{you see 1,40. Add 1,40 to 4,40,0, (and) you see 4,41,40. What is the square root?  You see 2,10.}\\
	(L8)\> \tabfill{Add 10, which was taken (in the above), to 2,10, (and) you see 2,20. What  should I put to  the (coefficient of provisional) area 28 which will give me 2,20?}\\
	(L9)\> \tabfill{Put down 5. Multiply 5 by 7, (and) you see 35. Subtract 5, which is subtracted  from the (provisional) length, from 35, (and)}\\
	(L10)\> \tabfill{you see 30. 30 is the (true) length. Multiply 5 (text; 5 of the length) by 4 of   the width, (and) you see 20. 20 is the (true) width (text; the length).} 
\end{tabbing}

\subsubsection*{Mathematical Calculations of First Problem}\label{SS-P1MC-SMT8} 
If we denote {\fontfamily{qpl}\selectfont u\v{s}}  ``the length'' by $x$ and {\fontfamily{qpl}\selectfont sag}  ``the width'' by $y$, the equations presented in line 1 would be as follows:

\begin{equation}\label{equ-8-a}
	\begin{dcases}
		xy  = 10,0  \\
		y + \dfrac{3}{4}y  = x + 5. 
	\end{dcases}
\end{equation}
Here the left-hand side of the second equation in  \cref{equ-8-a}  is the (provisional) width. In line 2, a transformation, namely 
\begin{equation}\label{equ-8-aa}
	y=4z,
\end{equation}
is introduced, and according to lines 2-3, the left-hand side of the second equation in  \cref{equ-8-a}  is changed into $7z$   giving   the following new equation 
\begin{equation}\label{equ-8-b}
	\begin{dcases}
		4xz  = 10,0  \\
		7z  = x + 5. 
	\end{dcases}
\end{equation}

To get the (true) length, according to line 4, the number 5 must be subtracted from $7z$, the (provisional) length. By the abbreviated steps  stated  in lines 4-6, the next calculations can be carried out  as follows:
\begin{align*}
	~~&~~4z(7z -5 )    = 10,0   \\ 
	\Longrightarrow~~&~~28z^2 - 20z  = 10,0  \\
	\Longrightarrow~~&~~(28z)^2 - 20\times(28z)  = (10,0)\times28\\
	\Longrightarrow~~&~~(28z)^2 - 20\times(28z)  = 4,40,0 
\end{align*} 
so we get  the new quadratic equation
\begin{equation}\label{equ-8-c}
	(28z)^2 - 20\times(28z)  =   4,40,0.
\end{equation}
As expected, this quadratic equation is solved in lines 6-8 by completing the square as follows:
\begin{align*}
	~~&~~(28z)^2 - 20\times(28z)  =   4,40,0\\
	\Longrightarrow~~&~~(28z)^2 - 2\times 10 \times(28z)  =   4,40,0\\
	\Longrightarrow~~&~~(28z)^2 - 2\times 10 \times(28z)+10^2  =   4,40,0+10^2\\
	\Longrightarrow~~&~~(28z -10 )^2  = 4,41,40 \\
	\Longrightarrow~~&~~28z- 10  = \sqrt{4,41,40}\\
	\Longrightarrow~~&~~28z-10 =2,10\\
	\Longrightarrow~~&~~28z = 2,20\\
	\Longrightarrow~~&~~z =\dfrac{1}{28}\times(2,20)\\
	\Longrightarrow~~&~~z =5.
\end{align*}  
Since the number 28 is an irregular number,   a special division is carried out in line 8, and it is   found that $z = 5$,  the positive solution of  \cref{equ-8-c}. This division  might have easily been done because of the fact that $2,20=140=7\times 4 \times 5$ is divided by $28=7\times 4$ and the scribe of this tablet would have not been worried about the irregular factor 7 in the numbers involved.  

Finally, in lines 9 and 10, the true length and width are obtained by using \cref{equ-8-aa} and \cref{equ-8-b} as follows:
\[ x = 7z - 5 = 35 - 5 = 30\]
and 
\[y = 4z = 4\times5 = 20, \]
and thus, the solutions of system of equations   \cref{equ-8-a}  are $x=30$ and $y=20$.

\subsubsection*{Second Problem}\label{SS-P2-SMT8} 
In the second problem, two technical terms peculiar to the Susa mathematical texts occur, that is, {\fontfamily{qpl}\selectfont gaba-4} and {\fontfamily{qpl}\selectfont \textit{t\={a}lukum}} (line 13). The former is a variation of {\fontfamily{qpl}\selectfont igi-4} ($=\frac{1}{4}$), probably because both {\fontfamily{qpl}\selectfont gaba} and {\fontfamily{qpl}\selectfont igi} are used as the logogram of {\fontfamily{qpl}\selectfont \textit{mihirtum}} ``counter part, front part''. The latter, which is derived from the verb {\fontfamily{qpl}\selectfont \textit{al\={a}kum}} ``to go; to multiply'', usually means ``gait, course'', but it is used as ``product (of two numbers)'' in our text.

The last three cuneiform signs of line 11 seem to be {\fontfamily{qpl}\selectfont \textit{i-\v{s}\'{i}-ma}} as in the \textbf{TMS}, but that  makes no sense. Judging from the mathematical content, we might expect them to be ({\fontfamily{qpl}\selectfont u\v{s} ugu sag}) {\fontfamily{qpl}\selectfont 5 dirig}, or the like, but   they do not fit  the signs. Therefore, on this occasion, we have given priority to the mathematical analysis of the text over the transliteration of the signs in our translation, while admitting that further  grammatical analysis might yield further insight as to their intended meaning.

\subsubsection*{Transliteration of Second Problem}\label{SSS-P2TI-SMT8} 
\begin{note1} 
	\underline{Obverse:  Lines 11-20} \\
	(L11)\hspace{1mm} [a-\v{s}\`{a} 10] 4-\textit{at} sag \textit{a-na} u\v{s}(sic) dah \textit{a-na} 1 \textit{a-li-ik a-na-di} [u]\v{s} ugu sag \textit{i-\v{s}\'{i}-ma}(?)\\
	(L12)\hspace{1mm} [za]-e 4 \textit{re-ba-ti ki-ma} sag gar \textit{re-ba-at} 4 \textit{le-q\'{e}} 1 \textit{ta-mar} 1 \textit{a-na} 1 \textit{a-li}-[\textit{ik}]\\
	(L13)\hspace{1mm} [1 \textit{ta-mar}] 4 gaba-4 gar 1 \textit{ta-lu-ka a-na} 4 dah 5 \textit{ta}-$<$\textit{mar}$>$ \textit{ki-ma} u\v{s} gar\\
	(L14)\hspace{1mm} [5 \textit{a-na}] \textit{wa-\c{s}\'{i}-ib} u\v{s} gar 5 u\v{s} \textit{a-na} 4 sag \textit{i-\v{s}\'{i}} 20 \textit{ta-mar} 20 a-\v{s}\`{a}\\
	(L15)\hspace{1mm} [\textit{a-na} 10] \textit{i}-[\textit{\v{s}\'{i}}] 3,20 \textit{ta-mar} 5 \textit{wa-\c{s}\'{i}-ib a-na} 4 sag \textit{i-\v{s}\'{i}} 20 \textit{ta-mar}\\
	(L16)\hspace{1mm} [1/2 \textit{he-pe} 10 \textit{ta-mar}] 10 nigin 1,40 \textit{ta-mar} 1,40 \textit{a-na} 3,20 dah 3,21,40 \textit{ta-mar}\\
	(L17)\hspace{1mm} [\'{i}b-si 3,21,40] \textit{le-q\'{e} i-na} 1,50 zi 1,40 \textit{ta-mar}\\
	(L18)\hspace{1mm} [igi-20 \textit{pu-\c{t}\'{u}-\'{u}r} 3 \textit{ta-mar} 3] \textit{a-na} 1,40 \textit{i-\v{s}\'{i}} 5 \textit{ta-m}[\textit{ar}] 5 \textit{a-na} 5 u\v{s}\\
	(L19)\hspace{1mm} [\textit{i-\v{s}\'{i}} 25 \textit{ta-mar} 5 \textit{wa-\c{s}\'{i}-ib} u\v{s} \textit{a}]-\textit{na} 25 dah 30 \textit{ta-mar} 30 u\v{s}\\
	(L20)\hspace{1mm} [5 \textit{a-na} 4 \textit{i-\v{s}\'{i}} 20 \textit{ta-mar}] 20 sag
\end{note1}

\subsubsection*{Translation of Second Problem}\label{SS-P2TR-SMT8} 
\underline{Obverse:  Lines 11-20}
\begin{tabbing}
	\hspace{12mm} \= \kill
	(L11)\> \tabfill{The area is 10,0. Add one fourth of the width to the width(text: length).  Multiply (one fourth of the width) by 1. I lay down that the length exceeds   the (provisional) width  by 5.}\\
	(L12)\> \tabfill{You, put down 4 of one fourth as the width. Take one fourth of 4, (and) you   see 1. Multiply 1 by 1, (and)}\\
	(L13)\> \tabfill{you see 1. Put down 4 of a fourth. Add the product 1 to 4, (and) you see 5.   Put down (this) as the length.} \\
	(L14)\> \tabfill{Put down 5 for that which is added for the length. Multiply the length 5 by   the width 4, (and) you see 20.} \\
	(L15)\> \tabfill{Multiply the (coefficient of provisional) area 20 by 10,0, (and) you see 3,20,0.   Multiply 5 , which is added, by the width 4, (and) you see 20.} \\
	(L16)\> \tabfill{Halve 20, (and) you see 10. Square 10, (and) you see 1,40. Add 1,40 to 3,20,0,  (and) you see 3,21,40.}\\
	(L17)\> \tabfill{Take the square root of 3,21,40, (and you see 1,50). Subtract (10) from 1,50,   (and) you see 1,40.} \\
	(L18)\> \tabfill{Make the reciprocal of 20, (and) you see 0;3. Multiply 0;3 by 1,40, (and) you   see 5.}\\
	(L19)\> \tabfill{Multiply 5 by the length 5, (and) you see 25. Add 5, which is added for the   length, to 25, (and) you see 30. 30 is the (true) length.}\\
	(L20)\> \tabfill{Multiply 5 by 4, (and) you see 20. 20 is the (true) width.}
\end{tabbing}

\subsubsection*{Mathematical Calculations of Second Problem}\label{SS-P2MC-SMT8} 
The mathematical structure of the second problem, which is similar to the first one, is as follows. Again, we denote {\fontfamily{qpl}\selectfont u\v{s}}  ``the length'' by $x$ and {\fontfamily{qpl}\selectfont sag}  ``the width'' by $y$. In this case,  line 11 gives us the following system of equations:
\begin{equation}\label{equ-8-d}
	\begin{dcases}
		xy = 10,0 \\
		y + \dfrac{1}{4}y + 5 = x. 
	\end{dcases}
\end{equation}
Lines 12-13  suggest  a new variable
\begin{equation}\label{equ-8-dd}
	y=4z,
\end{equation}
which gives $ y + \dfrac{1}{4}y = 5z$,
so  the first equation becomes 
\begin{align*}
	~~&~~xy =10,0\\
	\Longrightarrow~~&~~4z (5z + 5)  =10,0\\
	\Longrightarrow~~&~~    20z^2 + 20z  = 10,0.
\end{align*}
Thus, we get the following quadratic equation
\begin{equation}\label{equ-8-ddd}
	20z^2 + 20z  = 10,0.
\end{equation} 
According to lines 14-18,   the solution of   the previous quadratic equation by using completing the square:
\begin{align*}
	~~&~~20z^2 + 20z  = 10,0\\
	\Longrightarrow~~&~~(20z)^2 + 20 \times(20z)  =  20\times (10,0) \\
	\Longrightarrow~~&~~(20z)^2 + 2\times 10 \times(20z)   =3,20,0\\
	\Longrightarrow~~&~~(20z)^2 + 2\times 10 \times(20z)+ 10^2    =3,20,0+ 10^2 \\
	\Longrightarrow~~&~~(20z + 10 )^2  = 3,20,0 + 1,40  \\
	\Longrightarrow~~&~~(20z + 10 )^2  =  3,21,40\\
	\Longrightarrow~~&~~20z + 10    =  \sqrt{3,21,40}\\
	\Longrightarrow~~&~~20z + 10    =  1,50\\
	\Longrightarrow~~&~~20z    =  1,50- 10 \\
	\Longrightarrow~~&~~20z =1,40\\
	\Longrightarrow~~&~~z =\dfrac{1}{20}\times (1,40)\\
	\Longrightarrow~~&~~z =(0;3)\times (1,40)\\
	\Longrightarrow~~&~~z =5.
\end{align*}
According to lines 19-20,    the solutions of  \cref{equ-8-d} are obtained by using \cref{equ-8-dd}    as follows:
\[  x = 5z + 5 = 5\times 5+5=30\]
and 
\[y = 4z =4\times 5= 20.\]

\subsubsection*{Third Problem}\label{SS-P3-SMT8} 
Nothing  remains of the text regarding this problem  except for a few words:

\begin{tabbing}
	\hspace{12mm} \= \kill
	(L21)\> \tabfill{{\fontfamily{qpl}\selectfont $\cdots $ $\cdots $ $\cdots $    -\textit{ma} ul-gar 3 u\v{s} \textit{\`{u}} 4 sag}}\\
	\> \tabfill{``$\cdots $  and the sum of three times the length and four times the width.''}\index{length}\index{width}\\
	(L22)\> \tabfill{{\fontfamily{qpl}\selectfont $\cdots $ $\cdots $ $\cdots $ $\cdots $ \textit{ta-mar}}}\\
	\> \tabfill{``$\cdots $ you see $\cdots $''}
\end{tabbing}

\begin{remark}\label{rem-SMT8-a}
	The reader should be aware that other interpretations of \textbf{SMT No.\,8}\index{SMT No.h@\textbf{SMT No.\,8}} have been  given by H\o yrup  in  \cite{Hyp02} and Friberg in \cite{FA16}. 
\end{remark}

\subsection{SMT No.\,11}

\subsubsection*{First Problem}\label{SS-P1-SMT11} 

\subsubsection*{Transliteration of First Problem}\label{SSS-P1TI-SMT11} 

\begin{note1} 
	\underline{Obverse:  Lines 1-7}\\
	(L1)\hspace{2mm} ul-gar u\v{s} \textit{\`{u}} sag \textit{it-ti} u[\v{s} nigin]\\
	(L2)\hspace{2mm} \textit{i-na \v{s}\`{a} i-la-am} 2 zi \textit{lu-\'{u}} [a-\v{s}\`{a}]\\
	(L3)\hspace{2mm} za-e 20 a-r\'{a} \textit{\v{s}\`{a} e-li-ka \c{t}\`{a}}-[\textit{bu} gar]\\
	(L4)\hspace{2mm} igi-20 \textit{pu-\c{t}\'{u}-\'{u}r} 3 \textit{ta-mar} 3 \textit{a}-[\textit{na} 2]\\
	(L5)\hspace{2mm} \textit{\v{s}\`{a} ta-na-as-s\`{a}-hu i-\v{s}\'{i}} 6 \textit{ta-mar} [6 \textit{a-na} 20]\\
	(L6)\hspace{2mm} [\textit{\v{s}\`{a} iq}]-\textit{ri-bu} dah 26 \textit{ta-mar} 2[6 ul-gar]\\
	(L7)\hspace{2mm} [\textit{\v{s}\`{a} te}]-\textit{el-q\'{u}-\'{u}}
\end{note1}

\subsubsection*{Translation of First Problem}\label{SSS-P1TR-SMT11} 

(L1)\hspace{2mm} I multiplied the sum of a length and a width by the length.\\ 
(L2)\hspace{2mm} I subtracted 2,0 from the result which came out. Let it be the square of the length.\\ 
(L3)\hspace{2mm} You, put down the factor 20 which is convenient for you.\\
(L4)\hspace{2mm} Make the reciprocal of a number of 20, (and) you see 0;3.\\
(L5)\hspace{2mm} Multiply 0;3 by 2,0 which you subtract, (and) you see 6. \\
(L6)\hspace{2mm} Add 6 to 20 which was at (your) hand, (and) you see 26. \\
(L7)\hspace{2mm} 26 is the sum which you took.

\subsubsection*{Mathematical Calculations of First Problem}\label{SSS-P1MC-SMT11} 
If we denote {\fontfamily{qpl}\selectfont u\v{s}} ``the length'' by $x$ and {\fontfamily{qpl}\selectfont sag} ``the width'' by $y$, then the equation described in lines 1-2 is 
\begin{equation}\label{equ-SMT11-a}
	x(x + y) - 2,0 = x^2. 
\end{equation} 
This equation is obviously indeterminate. If  (according to line 3) we substitute $x = 20$ in  \cref{equ-SMT11-a}, we (according to lines 4-6) have
\begin{align*}
&~~ x(x+y) =2,0+x^2\\
\Longrightarrow~~&~~ x+y =\dfrac{2,0}{x}+x\\
\Longrightarrow~~&~~ 20+y = \dfrac{1}{20} \times(2,0)+20\\
\Longrightarrow~~&~~ 20+y =(0;3)\times(2,0)+20\\
\Longrightarrow~~&~~  20+y =26.
\end{align*}   
Although the subtraction $y = 26 - 20 = 6$ is expected in the remaining lines, the text suddenly ends at the beginning of line 7, leaving the blank space.

In line 3 the technical expression  {\fontfamily{qpl}\selectfont \textit{e-li-ka \c{t}\`{a}-bu}} ``(which) is convenient for you''    occurs, indicating that both $x$ and $y$ are indeterminate. Another example of this expression is   {\fontfamily{qpl}\selectfont u\v{s} sag-ki \textit{ma-la e-li-ia \c{t}\`{a}-bu \v{s}a-ka-nam} 1(\`{e}\v{s}e) iku \textit{ba-na-am}} ``To put down a length (and) a width such as are convenient for me, (and) to make the area  1 {\fontfamily{qpl}\selectfont \`{e}\v{s}e}, that is, $xy = 10,0$ ({\fontfamily{qpl}\selectfont sar})''\footnote{Note that {\fontfamily{qpl}\selectfont \`{e}\v{s}e}  and {\fontfamily{qpl}\selectfont sar} are area   units used in Babylonian mathematics such that 1 {\fontfamily{qpl}\selectfont   \`{e}\v{s}e} is 10,0 {\fontfamily{qpl}\selectfont sar} and 1 {\fontfamily{qpl}\selectfont sar} is almost equal to 36$m^2$.}. 
For a discussion on this  topic, see \cite{Mur07-1}.

\subsubsection*{Second Problem}\label{SS-P2-SMT11} 

\subsubsection*{Transliteration of Second Problem}\label{SSS-P2TI-SMT11} 

\begin{note1} 
\underline{Obverse, Bottom, Reverse:  Lines 8-26}\\
(L8)\hspace{2mm} [a-\v{s}\`{a} 20] a-r\'{a} \textit{\v{s}\`{a} e-li-i}[\textit{a \c{t}}]\textit{\`{a}-bu a-n}[\textit{a} u\v{s}]\\
(L9)\hspace{2mm}  [\textit{al-li-ik}] a-\v{s}\`{a} \textit{\`{u}} u\v{s} ul-ga[r] 46,[40 20 \textit{\v{s}\`{a}} u\v{s}]\\
(L10)\hspace{0mm}  [\textit{a-na} sag] dah 30 \textit{\v{s}\`{a}} sag \{\textit{\v{s}\`{a}} sag\}\\
(L11)\hspace{0mm}  [35 za]-e [3]5 \textit{a-na} 3[0 \textit{i-\v{s}\'{i}} 17,30 \textit{ta-mar}]\\

\underline{Bottom}\\
(L12)\hspace{0mm}  17,30 \textit{i-na} 46,40 zi 2[9,10 \textit{ta-mar}]\\
(L13)\hspace{0mm}  \textit{r}[\textit{e-i}]\textit{\v{s}-ka} \textit{li-ki-il tu-\'{u}r-m}[\textit{a}]\\

\underline{Reverse}\\
(L14)\hspace{0mm}  \textit{mi-na} dirig 20 dirig 20 dah gar 20 [nigin] \\
(L15)\hspace{0mm}  6,40 a-\v{s}\`{a} \textit{s\`{a}-ar-ru} 6,40 \textit{a}-[\textit{na} 29,10 \textit{i-\v{s}\'{i}}]\\
(L16)\hspace{0mm}  3,14,26,40 \textit{ta-mar} 30 \textit{a-na} [20 \textit{i-\v{s}\'{i}} 10 \textit{ta-mar}]\\
(L17)\hspace{0mm}  35 \textit{a-na} a-r\'{a} 20 \textit{i-\v{s}\'{i}} 11,40 \textit{t}[\textit{a-mar}]\\
(L18)\hspace{0mm}  \textit{\'{u}-ul} ul-gar 10 \textit{i-na} 11,40 zi [1,40 \textit{ta-mar} 1]\\
(L19)\hspace{0mm}  \textit{ka-aiia-ma}-[\textit{n}]\textit{a a-na} 1,[40] dah 1,1,40 [\textit{ta-mar}]\\
(L20)\hspace{0mm}  $ \frac{\text{1}}{\text{2}} $ 1,1,40 \textit{he-pe} 3[0],50 \textit{ta}-$<$\textit{mar}$>$ 30,[50 nigin 15,50,41,40 \textit{ta-mar}]\\
(L21)\hspace{0mm}  3,14,26,40 \textit{i}-[\textit{na}] 15,50,4[1,40 zi]\\
(L22)\hspace{0mm}  12,36,15 \textit{t}[\textit{a-ma}]\textit{r} \textit{mi-n}[\textit{a}] \'{i}[b-si 27,30 \'{i}b-si]\\
(L23)\hspace{0mm}  27,30 \textit{i-na} 30,[50] \textit{le-q\'{e}} zi [3,20 \textit{ta-mar}]\\
(L24)\hspace{0mm}  igi-6,40 \textit{pu-\c{t}\'{u}}-[\textit{\'{u}r}] 9 \textit{a-na} 3,[20 \textit{i-\v{s}\'{i}} 30 \textit{ta-mar}]\\
(L25)\hspace{0mm}  30 \textit{a-na} 1 \textit{i-\v{s}}[\textit{\'{i}}] 30 \textit{ta}-$<$\textit{mar}$>$ 30 [u\v{s}]\\
(L26)\hspace{0mm}  30 \textit{a-na} 40 \textit{i-\v{s}\'{i}} 20 \textit{ta-mar} [20 sag]
\end{note1}

\subsubsection*{Translation of Second Problem}\label{SSS-P2TR-SMT11} 

\begin{tabbing}
\hspace{13mm} \= \kill
(L8-9)\> \tabfill{The area ($xy$). I multiplied 0;20, the factor which is convenient for me, by the length. (Then) I added the area ($= 0;20xy$), the length, and 0;30 of the width, (and the result is) 0;46,40.}\\
(L10)\> \tabfill{I added 0;20 of the length to the width, (and the result is)}\\
(L11)\> \tabfill{0;35. You, multiply 0;35 by 0;30, (and) you see 0;17,30.}\\
(L12)\> \tabfill{Subtract 0;17,30 from 0;46,40, (and) you see 0;29,10.}\\
(L13)\> \tabfill{Let it hold your head (that is, memorize it). Return.}\\
(L14)\> \tabfill{What dose it (= 0;35) exceed? It exceeds 0;20 (in equation form: $y = 0;35 - 0;20x$). Put down 0;20 which was added (to the width). Square 0;20, (and)}\\
(L15)\> \tabfill{0;6,40, the provisional area. Multiply 0;6,40 by 0;29,10, (and)}\\
(L16)\> \tabfill{you see 0;3,14,26,40. Multiply 0;30 by 0;20, (and) you see 0;10.}\\
(L17)\> \tabfill{Multiply 0;35 by the factor 0;20, (and) you see 0;11,40.} \\
(L18)\> \tabfill{Do not add! Subtract 0;10 from 0;11,40, (and) you see 0;1,40.} \\
(L19)\> \tabfill{Add the ``regular (number)'' 1 to 0;1,40, (and) you see 1;1,40.}\\
(L20)\> \tabfill{Halve 1;1,40, (and) you see 0;30,50. Square 0;30,50, (and) you see 0;15,50,41,40.}\\
(L21)\> \tabfill{Subtract 0;3,14,26,40 from 0;15,50,41,40, (and)}\\
(L22)\> \tabfill{you see 0;12,36,15. What is the square root (of it)? 0;27,30 is the square root.}\\
(L23)\> \tabfill{Subtract 0;27,30 from 0;30,50, which was taken (from 1;1,40 by halving), (and) you see 0;3,20.}\\
(L24)\> \tabfill{Make the reciprocal of 0;6,40. Multiply (the result) 9 by 0;3,20, (and) you see 0;30.} \\
(L25)\> \tabfill{Multiply 0;30 by 1, (and) you see 0;30. 0;30 is the length.} \\
(L26)\> \tabfill{Multiply 0;30 by 0;40, (and) you see 0;20. 0;20 is the width.}
\end{tabbing}

\subsubsection*{Mathematical Calculations of Second Problem}\label{SSS-P2MC-SMT11}  
The statement of the second problem in lines 8-10 is apparently in disorder, that is,
{\fontfamily{qpl}\selectfont 30 \textit{\v{s}\`{a}} sag} (``0;30 of the width'')  in line 10, of which {\fontfamily{qpl}\selectfont \textit{\v{s}\`{a}} sag} erroneously written down
twice, should be in line 9:

\noindent
(L9) {\fontfamily{qpl}\selectfont $\cdots$ a-\v{s}\`{a} u\v{s} \textit{\`{u}} 30 \textit{\v{s}\`{a}} sag ul-gar} \\
``$\cdots$ I added the area, the length, and 0;30 of the width''

\noindent
Moreover, it is strange that the phrase {\fontfamily{qpl}\selectfont a-r\'{a} \textit{\v{s}\`{a}} \textit{e-li-ia \c{t}\`{a}-bu}} ``the factor which is convenient for me'' in line 8, does not play any role in the problem. It is obvious, however, that the second problem treats the following system of simultaneous linear equations:
\begin{equation}\label{equ-SMT11-b}
\begin{dcases}
	(0;20)xy + x + (0;30)y = 0;46,40\\
	(0;20)x + y = 0;35
\end{dcases}
\end{equation}
By substituting  $y=0;35-(0;20)x$ in the first equation of \cref{equ-SMT11-b}, we get
\begin{align*}
(0;20)x\times \Big( 0;35-(0;20)x\Big)+ x  + (0;30)\times\Big(0;35-(0;20)x \Big) = 0;46,40 
\end{align*}
or 
\begin{align*}
\Big( (0;20)\times (0;35)\Big)x-(0;20)^2x^2  + x  +   (0;30)\times (0;35) -\Big( (0;20)\times (0;30)\Big) x= 0;46,40
\end{align*}
which gives us the following quadratic equation
\begin{equation}\label{equ-SMT11-c}
\begin{split}
	(0;20)^2x^2   -     \Big ( (0;20)\times (0;35)  -(0;20)&\times (0;30)+1 \Big ) x\\
	&+    0;46,40-(0;30)\times (0;35) =0.
\end{split}
\end{equation}

\noindent  
The following successive calculations are given  in lines 11-19
\begin{align*}
(\mathrm{L}11)&: \ \    (0;35)\times(0;30) = 0;17,30 \\
(\mathrm{L}12)&: \ \  0;46,40 - 0;17,30 = 0;29,10 \\
(\mathrm{L}14,15)&: \ \  (0;20)^2 = 0;6,40 \\
(\mathrm{L}15,16)&: \ \ (0;6,40)\times(0;29,10) = 0;3,14,26,40 \\
(\mathrm{L}16)&: \ \ (0;30)\times(0;20) = 0;10 \\
(\mathrm{L}17)&:  \ \ (0;35)\times(0;20) = 0;11,40 \\
(\mathrm{L}18)&:  \ \ 0;11,40 - 0;10 = 0;1,40 \\
(\mathrm{L}19)&:  \ \ 0;1,40 + 1 = 1;1,40. 
\end{align*}
These suggest the simplification of equation \cref{equ-SMT11-c} as follows:
\begin{align*}
~~&~~(0;6,40)x^2  -      ( 1;1,40 )  x +    0;29,10 =0\\
\Longrightarrow~~&~~(0;6,40)^2x^2  -       \Big((0;6,40)\times ( 1;1,40 )\Big)  x +     (0;6,40)\times(0;29,10) =0 
\end{align*}
which gives the following new quadratic equation
\begin{equation}\label{equ-SMT11-d}
\Big( (0;6,40)x\Big)^2 - (1;1,40)\times \Big((0;6,40)x\Big)   + 0;3,14,26,40 = 0.
\end{equation} 
Note that  the number 0;6,40 is called ``(of) the provisional area''.

Next, (according to lines 20-24)  the usual method of Babylonian mathematics  to solve equation \cref{equ-SMT11-d},  that is,  completing the square,  starts:
\begin{align*}
	~~&~~\Big( (0;6,40)x\Big)^2 - (1;1,40) \times  \Big((0;6,40)x\Big)    + 0;3,14,26,40 = 0\\
	\Longrightarrow~~&~~\Big( (0;6,40)x\Big)^2 - 2\times (0;30,50) \times \Big((0;6,40)x\Big)   + 0;3,14,26,40 = 0\\
	\Longrightarrow~~&~~\Big( (0;6,40)x\Big)^2 - 2\times (0;30,50)\times \Big((0;6,40)x\Big)    + 0;3,14,26,40\\
	~~&~~\qquad \qquad\qquad\qquad\qquad\qquad\qquad\qquad~~~~+(0;30,50)^2= (0;30,50)^2\\
	\Longrightarrow~~&~~\Big( (0;6,40)x\Big)^2 - 2\times (0;30,50)\times \Big((0;6,40)x\Big)    + (0;30,50)^2 \\
	~~&~~\qquad \qquad\qquad\qquad\qquad\qquad\qquad\qquad= (0;30,50)^2-0;3,14,26,40\\
	\Longrightarrow~~&~~\Big( (0;6,40)x\Big)^2 - 2\times (0;30,50)\times \Big((0;6,40)x\Big)    +  (0;30,50)^2\\
	~~&~~\qquad \qquad\qquad\qquad\qquad\qquad\qquad~= 0;15,50,41,40-0;3,14,26,40\\
	\Longrightarrow~~&~~\Big( (0;6,40)x  -  0;30,50 \Big)^2  = 0;12,36,15. 
\end{align*}

\noindent 
From here on, it is implicitly supposed that 0;30,50 is greater than the value of $(0;6,40)x$,  because   Susa scribes did not know the concept of negative numbers and    they always subtracted  a smaller number from a greater one.    Thus  we should switch these two  terms in the left-hand side of the last equation and continue as follows:
\begin{align*}
	~~&~~\Big(   0;30,50 - (0;6,40)x  \Big)^2   = 0;12,36,15\\
	\Longrightarrow~~&~~0;30,50 - (0;6,40)x     = \sqrt{0;12,36,15}\\
	\Longrightarrow~~&~~0;30,50 - (0;6,40)x     = \sqrt{(0;27,30)^2}\\
	\Longrightarrow~~&~~(0;6,40)x  = 0;30,50 - 0;27,30\\
	\Longrightarrow~~&~~(0;6,40)x  = 0;3,20\\
	\Longrightarrow~~&~~x =\dfrac{1}{(0;6,40)}\times (0;3,20)\\
	\Longrightarrow~~&~~x =9\times (0;3,20)\\
	\Longrightarrow~~&~~x =0;30.
\end{align*}

In last two  lines 25-26, some confusing calculations are carried out by the scribe as follows:
\[  \text{the length\index{length}}  =(0;30)\times 1 = 0;30   \]
and
\[   \text{the width\index{width}}  =(0;30)\times(0;40) = 0;20.\]
In fact, the former gives redundant  information and the latter is wrong and should be  replaced with the following 
\[  y = 0;35 - (0;20)x = 0;35-(0;20)\times(0;30)=0;35-0;10=0;25.\] 
Therefore, as   expected, the solutions of system of linear equations \cref{equ-SMT11-b} are $x=0;30$ and $y=0;20$.

\begin{remark}\label{rem-SMT11-a}
	In line 23, the verbal adjective {\fontfamily{qpl}\selectfont \textit{leq\^{u}m}} ``to be taken''   is used as    {\fontfamily{qpl}\selectfont 30,50 \textit{le-q\'{i}}} ``0;30,50 which was taken (from 1;1,40 by halving)''.  
	We should not confuse this technical expression with {\fontfamily{qpl}\selectfont \textit{le-q\'{e}}} ``take!'', which often occurs in mathematical texts. We saw the same term   in \textbf{SMT No.\,8}, line 8,  in the same mathematical situation (see \cref{SS-P1-SMT8}).
\end{remark}

\subsection{SMT No.\,17}

\subsubsection*{Transliteration}\label{SS-TI-SMT17}

\begin{note1}
	\underline{Obverse:  Lines 1-12}\\
	(L1)\hspace{2mm} u\v{s} ki sag \textit{u\v{s}-ta-ki-il-ma} a-\v{s}\`{a} \textit{ab-ni}\\
	(L2)\hspace{2mm} \textit{\v{s}\`{a}} u\v{s} ugu sag dirig \textit{u\v{s}-ta-ki-il-ma}\\
	(L3)\hspace{2mm} \textit{a-na} a-\v{s}\`{a} dah-\textit{ma} 20 a-\v{s}\`{a} 4-\textit{at} u\v{s}\\
	(L4)\hspace{2mm} \textit{\`{u}} 4-\textit{at} sag ul-gar-\textit{ma} 15 u\v{s} \textit{\`{u}} sag \textit{mi-nu-um}\\
	(L5)\hspace{2mm} [z]a-e \textit{a\v{s}-\v{s}um} 4\textit{-tum qa-bi-a-ku-u}m 15 ul-gar\\
	(L6)\hspace{2mm} \textit{a-di} 4 t\'{u}m-\textit{ma} 1 \textit{ta-mar} 1 ul-gar u\v{s} \textit{\`{u}} sag\\ 
	(L7)\hspace{2mm} 1 ul-gar \textit{\v{s}u-ta-ki-il-ma} 1 \textit{ta-mar}\\
	(L8)\hspace{2mm} 20 a-\v{s}\`{a} \textit{i-na} 1 zi-\textit{ma} 40 \textit{ta-mar tu-úr-ma}\\
	(L9)\hspace{2mm} [40] \textit{ki-ma} 1 u\v{s} 1 \textit{ki-ma} $<$1$>$ sag gar\\
	(L10)\hspace{0mm} [40] \textit{ta-mar a\v{s}-\v{s}um} dirig \textit{\v{s}u-ta-ku-lu qa-bu-ku-um}\\
	(L11)\hspace{0mm} [40] \textit{ù} 1 \textit{\v{s}\`{a} ki-ma} u\v{s} \textit{\`{u}} sag gar ul-gar-\textit{ma}\\
	(L12)\hspace{0mm}  [\textit{i-na} 2 u\v{s} zi-\textit{ma} dirig \textit{\`{u}} dirig] \textit{i-di \v{s}u-ta-ki-il-ma}\\
	
	\underline{Reverse:  Lines 1-9}\\
	(L1)\hspace{2mm} [40 \textit{ki-ma} u\v{s}] \textit{a-na} 1 sag \textit{i-}[\textit{\v{s}i}]-\textit{ma} 40 \textit{ta-mar} $<$40 u\v{s}$>$\\
	(L2)\hspace{2mm} [\textit{i-na} ul-ga]r zi-\textit{ma} 20 \textit{ta-mar} 20 sag 40 a-r\'{a}
	\\
	(L3)\hspace{2mm} [\textit{a-na}] 2 \textit{\v{s}\`{a} ki-ma} dirig \textit{i-\v{s}i-ma} 1,20 \textit{ta-mar}\\  
	(L4)\hspace{2mm} 1 ul-gar \textit{i-na} 1,20 zi-\textit{ma} 20 \textit{ta-mar} 20 dirig\\
	(L5)\hspace{2mm} 40 u\v{s} \textit{a-na} 20 sag \textit{i-\v{s}i-ma} 13,20 \textit{ta-mar} 13,20 a-\v{s}\`{a} \textit{ke-nu}\\
	(L6)\hspace{2mm} 20 dirig \textit{\v{s}u-ta-ki-il-ma} 6,40 \textit{ta-mar}\\
	(L7)\hspace{2mm} 6,40 \textit{a-na} 13,20 a-\v{s}\`{a} dah-\textit{ma} 20 \textit{ta-mar}\\
	(L8)\hspace{2mm} 20 ul-gar a-\v{s}\`{a} \textit{\`{u}} a-\v{s}\`{a} dirig
	\\
	(L9)\hspace{2mm} \textit{ki-a-am ma-ak-\c{s}a-ar-\v{s}u}
	
\end{note1}

\subsubsection*{Translation}\label{SS-TR-SMT17}  
\underline{Obverse:  Lines 1-12}
\begin{tabbing}
	\hspace{14mm} \= \kill 
	(L1)\> \tabfill{I multiplied the length and the width together, and I made the area.}\\
	(L2)\> \tabfill{I squared that by which the length exceeds the width, and} \\ 
	(L3)\> \tabfill{I added (it) to the area, and (the result is) 0;20.}\\
	(L4)\> \tabfill{I added one fourth of the length and one fourth of the width together, and (the result is) 0;15. What are the length and the width?}\\
	(L5)\> \tabfill{You, since one fourth is said to you,}\\
	(L6)\> \tabfill{multiply 0;15 of the sum by 4, and you see 1. 1 is the sum of the length and the width.}\\
	(L7)\> \tabfill{Square 1 of the sum, and you see 1.}\\ 
	(L8)\> \tabfill{Subtract 0;20 of the areas from 1, and you see 0;40. Return.}\\
	(L9)\> \tabfill{Put down 0;40 as one length (and) 1 as (one) width. Multiply (them) together, and}\\
	(L10)\> \tabfill{you see 0;40. Since squaring the excess is said to you,}\\
	(L11)\> \tabfill{put down 0;40 and 1 as (one) length and (one) width. The sum (of the length and the width),}\\
	(L12)\> \tabfill{subtract (it) from two-length, and write down the excess twice, (and then) square (the excess), and}
\end{tabbing}

\noindent
\underline{Reverse:  Lines 1-9}
\begin{tabbing}
	\hspace{14mm} \= \kill 
	(L1)\> \tabfill{multiply 0;40 as the length by 1 of (one) width, and you see 0;40. 0;40 is the length.}\\
	(L2)\> \tabfill{Subtract (it) from the sum (of the length and the width), and you see 0;20. 0;20 is the width. The factor 0;40 (that is 0;40 of one length),}\\
	(L3)\> \tabfill{multiply (it) by 2 that is (used in calculation of) the excess, and you see 1;20.}\\ 
	(L4)\> \tabfill{Subtract 1 of the sum from 1;20, and you see 0;20. 0;20 is the excess.}\\ 
	(L5)\> \tabfill{Multiply 0;40 of the length by 0;20 of the width, and you see 0;13,20. 0;13,20 is the true area.}\\
	(L6)\> \tabfill{Square 0;20 of the excess, and you see 0;6,40.}\\ 
	(L7)\> \tabfill{Add 0;6,40 to 0;13,20, and you see 0;20.}\\
	(L8)\> \tabfill{0;20 is the sum of the area and the square of the excess.}\\
	(L9)\> \tabfill{Such is its factorization (method).}
\end{tabbing}

\subsubsection*{Mathematical Calculations}\label{SS-MC-SMT17} 
Let $x$ and $y$ denote the length and  the width respectively. Lines 1-4 give us the following system of equations:

\begin{equation}\label{equ-SMT17-a}
	\begin{dcases}
		xy+(x-y)^2=0;20\\
		\dfrac{1}{4} x+\dfrac{1}{4}y=0;15.
	\end{dcases}
\end{equation} 

\noindent
After readying the translation more carefully, we can see that  the scribe has used the following algebraic identity to solve the previous system of equations:
\begin{equation}\label{equ-SMT17-b}
	(x+y)^2-\left[ xy+(x-y)^2 \right]=3xy.
\end{equation} 

\noindent
In fact, (according to lines 5-9) we first multiply both sides of the second equation in \cref{equ-SMT17-a} to get
\begin{equation}\label{equ-SMT17-c}
	x+y =4\times (0;15)=1
\end{equation} 

\noindent
which immediately implies that
\begin{equation}\label{equ-SMT17-d}
	(x+y)^2 =1.
\end{equation}

\noindent
Now, it follows from \cref{equ-SMT17-a}, \cref{equ-SMT17-b} and \cref{equ-SMT17-d} that
\begin{align*}
	&~~  3xy = (x+y)^2-\left[ xy+(x-y)^2 \right] \\
	\Longrightarrow~~&~~  3xy=  1-0;20 \\
	\Longrightarrow~~&~~    3xy =  0;40
\end{align*}
which implies that
\begin{equation}\label{equ-SMT17-e}
	3xy =  0;40.
\end{equation} 

From here on, we may expect the scribe to find the value of $xy$ and then apply the usual Babylonian method to find $x$ and $y$ (see \cref{rem-SMT17-aa} below). Unexpectedly, after lines 10-12,  the scribe  changes the direction of the  common  solution completely. He appears to assume that
\[x=0;40\ \ \  \text{and}\ \ \  3y=1.\]
which clearly comes from  the equality $x\times(3y)=3xy=0;40=(0;40)\times 1$.  Next, in lines 1-4 on the reverse, the scribe does the following calculations:
\[ x=(0;40)\times 1=0;40 \]
and
\[ y=(x+y)-x=1-0;40=0;20 \]
which gives the value of $x,y$. 
Then, he calculates the value of $x-y$ as follows:
\[ x-y=2x-(x+y)=2\times(0;40)-1=1;20-1=0;20. \]
Finally, in lines 5-8 on the reverse, he uses the previous values to verify that the numbers $x=0;40$ and $y=0;20$ satisfy the system of equations \cref{equ-SMT17-a}:
\begin{align*}
	xy+(x-y)^2&=(0;40)\times (0;20)+(0;40-0;20)^2\\	
	&=0;13,20+(0;20)^2\\
	&=0;13,20+0;6,40\\
	&=0;20.
\end{align*}

\begin{remark}\label{rem-SMT17-aa}
	We can use the usual method   to find $x$ and $y$, as it follows from \cref{equ-SMT17-e} that
	\begin{align*}
		&~~  3xy = 0;40\\
		\Longrightarrow~~&~~   xy=  \frac{1}{3}\times (0;40) \\
		\Longrightarrow~~&~~     xy =(0;20)\times  0;40
	\end{align*}
	or
	\begin{equation}\label{equ-SMT17-f}
		xy =0;13,20
	\end{equation}  
	and then by \cref{equ-SMT17-d} and \cref{equ-SMT17-f} we can write
	\begin{align*}
		\frac{(x-y)}{2} &=\sqrt{\left(\frac{x+y}{2}\right)^2-xy}  \\
		&  =  \sqrt{\left(\frac{1}{2}\right)^2-0;13,20} \\
		&   =  \sqrt{0;15-0;13,20} \\
		&    =   \sqrt{0;1,40} \\
		&   =  0;10
	\end{align*}
	
	\noindent
	so
	\begin{equation}\label{equ-SMT17-g}
		\frac{x-y}{2} =0;10.
	\end{equation} 
	
	\noindent
	By \cref{equ-SMT17-c} and \cref{equ-SMT17-g}, we finally get
	\[ x=\frac{x+y}{2}+\frac{x-y}{2} =0;30+0;10=0;40  \]
	and
	\[ y=\frac{x+y}{2}-\frac{x-y}{2} =0;30-0;10=0;20.  \]
\end{remark}

\begin{remark}\label{rem-SMT17-a}
	Note that in \cref{equ-SMT17-b}, we are actually using  the well-known algebraic identity
	$$(A+B)^2-(A-B)^2=4AB. $$
	The factorization of algebraic expressions such as \cref{equ-SMT17-b} was one of the methods Babylonian used to solve   algebraic equations. For more complicated examples of this method, see \cite{Mur19}. 
\end{remark}

\subsection{SMT No.\,19}

\subsubsection*{Transliteration}\label{SSS-P2TI-SMT19}

\begin{note1}
	\underline{Reverse:  Lines 1-13}\\
	(L1)\hspace{2mm} bar-d\'{a} 20 a-\v{s}à \{\v{s}à\} u\v{s} nigin a-\v{s}\`{a} u\v{s} \textit{i-la-am-ma} ki bar-d\'{a} nigin-\textit{ma}
	\\
	(L2)\hspace{2mm} 14,48,53,20 u\v{s} sag \textit{\`{u}} bar-d\'{a} \textit{mi-nu}
	\\
	(L3)\hspace{2mm} za-e 20 a-\v{s}\`{a} ni[gin] 6,40 \textit{ta-mar re-i\v{s}-ka li-ki-il tu-\'{u}r-ma}
	\\
	(L4)\hspace{2mm} 14,48,5[3,20 nig]in 3,39,[28],44(sic),$<$27,24$>$,26,40 \textit{ta-mar}\\
	(L5)\hspace{2mm} 1/2 6,40 \textit{\v{s}\`{a} re-i\v{s}}-$<$\textit{ka}$>$ \textit{\'{u}-ki}-[\textit{lu h}]\textit{e-pe} [3],20 \textit{ta-mar}
	\\
	(L6)\hspace{2mm} 3,20 nigin 11,6,40 \textit{ta-mar} 11,[6,40]
	\\
	(L7)\hspace{2mm} \textit{a-na} 3,39,28,44(sic),$<$27,24$>$,26,[40] dah-[\textit{ma}]
	\\
	(L8)\hspace{2mm} 3,50,36(sic),43(sic),34(sic),26,40 \textit{ta-mar mi-na} \'{i}[b-s]i
	\\
	(L9)\hspace{2mm} 15,11,6,40 íb-si \textit{i-na} 15,11,6,40\\
	(L10)\hspace{0mm} 3,20 \textit{ta-ki-il-ta} 2 zi 11,51,6,40 \textit{ta-mar} \textit{mi-na} íb-si\\
	(L11)\hspace{0mm} 26,40 \'{i}b-si igi-26,40 \textit{pu-\c{t}\'{u}-\'{u}r} 2,15 \textit{ta-mar}
	\\
	(L12)\hspace{0mm} 2,15 \textit{a-na} 6,40 \textit{ta-ki-il-ti \v{s}\`{a} re-i\v{s}-ka u-ki-lu i-\v{s}\'{i}}
	\\
	(L13)\hspace{0mm} [1]5 \textit{ta-mar mi-na} \'{i}b-si 30 \'{i}b-si 30 sag
	\\
	
\end{note1}

\subsubsection*{Translation}\label{SSS-P2TR-SMT19}

\underline{Reverse:  Lines 1-13}
\begin{tabbing}
	\hspace{14mm} \= \kill 
	(L1)\> \tabfill{The diagonal. 20,0 is the area. I squared the length, (and) the (squared) area (multiplied by) the length came up, and I multiplied (the product) and the diagonal together, and (the result is)}\\
	(L2)\> \tabfill{14,48,53,20. What are the length, the width, and the diagonal?}\\
	(L3)\> \tabfill{You, square the area 20,0, (and) you see 6,40,0,0. Let it hold your head (memorize it). Return.}\\ 
	(L4)\> \tabfill{Square 14,48,53,20, (and) you see 3,39,28,43,27,24,26,40.}\\
	(L5)\> \tabfill{Halve 6,40,0,0 that held your head, (and) you see 3,20,0,0.}\\
	(L6)\> \tabfill{Square 3,20,0,0, (and) you see 11,6,40,0,0,0,0.}\\
	(L7)\> \tabfill{Add (it) to 3,39,28,43,27,24,26,40, and}\\
	(L8)\> \tabfill{you see 3,50,35,23,27,24,26,40. What is the square root?}\\
	(L9)\> \tabfill{15,11,6,40 is the square root. From 15,11,6,40,}\\
	(L10)\> \tabfill{subtract 3,20,0,0 that was (used in) the second step of completing the square, (and) you see 11,51,6,40. What is the square root?}\index{square root}\\
	(L11)\> \tabfill{26,40 is the square root. Make the reciprocal of 26,40, (and) you see 0;0,2,15.} \\ 
	(L12)\> \tabfill{Multiply 0;0,2,15 by 6,40,0,0 that was (used in) completing the square (and) held your head, (and)}\\ 
	(L13)\> \tabfill{you see 15,0. What is the square root? 30 is the square root. 30 is the width.}
\end{tabbing}

\subsubsection*{Mathematical Calculations}\label{SSS-P2MC-SMT19}   
There are three quantities in this problem: the length, the width, and the diagonal of a rectangle   which we denote by   $x$, $y$, and $d$  respectively. Note that by the Pythagorean theorem, $d=\sqrt{x^2+y^2}$.  Lines 1-2 give the following system of equations:
\begin{equation}\label{equ-SMT19-d}
	\begin{dcases}
		xy=20,0\\
		x^3\sqrt{x^2+y^2}=14,48,53,20.
	\end{dcases}
\end{equation}

The approach taken by the scribe here to solve this system of equations is creative. One of the standard methods   is to apply the elimination method  and to find the value of $x$ or $y$  with respect to one another from the first equation and then substitute it in the second one. This gives rise to an equation of degree 8 with respect to $x$ or $y$  which fortunately can be transformed into a quadratic equation by a proper transformation and  solved by completing the square. If we set, for example,  $y=(20,0)/x$ in the second equation and simplify, we get the following octic equation:
\begin{equation}\label{equ-SMT19-da}
	x^8+(6,40,0,0)x^4=3,39,28,43,27,24,26,40	
\end{equation}
which changes to a standard quadratic equation if we use the transformation  $z=x^4$: 
\begin{equation}\label{equ-SMT19-db}
	z^2+(6,40,0,0)z=3,39,28,43,27,24,26,40
\end{equation}
By completing the square, one  can  solve this quadratic equation and find the value of $z$ satisfying \cref{equ-SMT19-db}. By taking the fourth root of the obtained number, we get the solution of \cref{equ-SMT19-da}.

Interestingly enough,   the scribe of this tablet seems to have adopted another approach to solve this system of equations which avoids eliminating any variables.   Close analysis of the text reveals that the scribe  has computed  the value of the following algebraic expression:
$$ x^4+\dfrac{x^2y^2}{2}. $$
To do so, as is known from the text, he apparently has  computed the value of the square of this algebraic expression by using the two equations in \cref{equ-SMT19-d}. First, note that
\begin{align*}
	\left(x^4+\dfrac{x^2y^2}{2}\right)^2&=(x^4)^2+2(x^4)\left(\dfrac{x^2y^2}{2}\right)+\left(\dfrac{x^2y^2}{2}\right)^2 \\
	&= x^8+x^6y^2+\dfrac{x^4y^4}{4} \\
	&= x^6(x^2+y^2)+\dfrac{x^4y^4}{4}
\end{align*} 
so
\begin{equation}\label{equ-SMT19-e}
	\left(x^4+\dfrac{x^2y^2}{2}\right)^2= x^6(x^2+y^2)+\dfrac{x^4y^4}{4}.
\end{equation} 
To find the value of the left-hand side of \cref{equ-SMT19-e}, the scribe tries to compute  the value   of each term  in the right-hand side of \cref{equ-SMT19-e}. This is what he has done  in lines 3-8,   according to which  he  squares both sides of equations in \cref{equ-SMT19-d} to get
\[(xy)^2=(20,0)^2=6,40,0,0\]
and
\[\left(x^3\sqrt{x^2+y^2}\right)^2=(14,48,53,20)^2=3,39,28,43,27,24,26,40.\] 
These two give us the following equations:
\begin{equation}\label{equ-SMT19-f}
	\begin{dcases}
		x^2y^2= 6,40,0,0\\
		x^6(x^2+y^2)= 3,39,28,43,27,24,26,40. 
	\end{dcases}
\end{equation}
Next,  according to  lines 5-6, he   halves   both sides of the first equation in \cref{equ-SMT19-f} to obtain
\begin{align*}
	&~~  x^2y^2= 6,40,0,0 \\
	\Longrightarrow~~&~~    \dfrac{x^2y^2}{2}= \dfrac{1}{2}\times (6,40,0,0)\\
	\Longrightarrow~~&~~     \dfrac{x^2y^2}{2}=3,20,0,0
\end{align*} 
so
\begin{equation}\label{equ-SMT19-g}
	\dfrac{x^2y^2}{2} =3,20,0,0.
\end{equation}
Then, he squares both sides of \cref{equ-SMT19-g} to get
\begin{align*}
	&~~  \dfrac{x^2y^2}{2}=3,20,0,0 \\
	\Longrightarrow~~&~~    \left(\dfrac{x^2y^2}{2}\right)^2= (3,20,0,0)^{2}\\
	\Longrightarrow~~&~~      \left(\dfrac{x^2y^2}{2}\right)^2=11,6,40,0,0,0,0,
\end{align*} 
so
\begin{equation}\label{equ-SMT19-h}
\dfrac{x^4y^4}{4} =11,6,40,0,0,0,0.
\end{equation}
Now,  according to lines 7-8, it follows from  \cref{equ-SMT19-e}, \cref{equ-SMT19-f}  and \cref{equ-SMT19-h}  that\footnote{It seems that the scribe has used the prime factorization $3,50,35,23,27,24,26,40= 2^{14}\times 5^{8}\times 41^2$ to find the square root of this huge number.}
\begin{align*}
	x^4+\dfrac{x^2y^2}{2}&=  \sqrt{\left(x^4+\dfrac{x^2y^2}{2}\right)^2} \\
	&= \sqrt{x^6(x^2+y^2)+\dfrac{x^4y^4}{4}} \\
	& = \sqrt{3,39,28,43,27,24,26,40+11,6,40,0,0,0,0}  \\
	& =\sqrt{3,50,35,23,27,24,26,40}  \\
	& =\sqrt{(15,11,6,40)^2}  \\
	&  =15,11,6,40.
\end{align*}  
Thus we get
\begin{equation}\label{equ-SMT19-i}
	x^4+\dfrac{x^2y^2}{2}=  15,11,6,40.
\end{equation} 
According to lines 9-10,  it follows from \cref{equ-SMT19-g} and  \cref{equ-SMT19-i} that
\begin{align*}
	&~~  x^4+\dfrac{x^2y^2}{2}=  15,11,6,40 \\
	\Longrightarrow~~&~~     x^4=15,11,6,40-\dfrac{x^2y^2}{2}  \\
	\Longrightarrow~~&~~     x^4=15,11,6,40-3,20,0,0\\
	\Longrightarrow~~&~~     x^4=11,51,6,40\\
	\Longrightarrow~~&~~     x^2=\sqrt{11,51,6,40}\\
	\Longrightarrow~~&~~     x^2=\sqrt{(26,40)^2}\\
	\Longrightarrow~~&~~     x^2=26,40.
\end{align*} 
Hence
\begin{equation}\label{equ-SMT19-j}
	x^2=  26,40.
\end{equation} 
According to lines 11-13, we can find $y$ by using  \cref{equ-SMT19-f} and  \cref{equ-SMT19-j}  as follows:
\begin{align*}
	&~~  x^2y^2=  6,40,0,0 \\
	\Longrightarrow~~&~~   y^2= \dfrac{1}{x^2}\times (6,40,0,0) \\
	\Longrightarrow~~&~~     y^2=\dfrac{1}{(26,40)} \times (6,40,0,0) \\
	\Longrightarrow~~&~~       y^2=(0;0,2,15)\times (6,40,0,0)\\
	\Longrightarrow~~&~~       y^2=15,0\\
	\Longrightarrow~~&~~       y =\sqrt{15,0}
\end{align*}
so
\begin{equation}\label{equ-SMT19-k}
	y =  30.
\end{equation} 
The scribe has apparently forgotten to compute the value of $x$ which easily can be obtained by taking the square root of \cref{equ-SMT19-j}:
\begin{equation}\label{equ-SMT19-l}
	x=\sqrt{26,40}=40. 
\end{equation}  
Also note that the value of $d$ is easily obtained from \cref{equ-SMT19-k} and \cref{equ-SMT19-l}:
\[ d= \sqrt{x^2+y^2}=\sqrt{26,40+15,0}=\sqrt{41,40}=50 \] 
therefore, the Pythagorean triple treated in this problem is $(30,40,50)=(10\times 3, 10\times 4,10\times 5) $.

\begin{remark}\label{rem-SMT19-a}
The reader should note that another interpretation of this text has been given by H\o yrup in \cite{Hyp02} (see also \cite{Fri07-2}).
\end{remark}

\section{Conclusion} 
That Susa scribes were familiar with systems of equations and how to solve such equations is clear from our analysis. They    used different techniques such as the change of variables, substitution, and algebraic identities, which demonstrate   the  mathematical skill employed by the Susa scribes in dealing with systems of equations. 

 Among the texts we have considered here, \textbf{SMT No.\,19} is of particular interest. Two rare phenomena occur in this text.  Firstly, the   complicated system  of algebraic equations treated here leads  to an unusual  algebraic equation  of degree eight, which we have not seen in any other Babylonian mathematical texts.  Secondly,  the scribe's remarkable ability   to manipulate   huge numbers, particularly  in the second problem,  sheds additional light on  the range of  mathematical skills available to at least  some of the Susa scribes. They could carry out  arithmetical computations involving very large numbers. For example, finding the square root of a significant sexagesimal  number such as $3,50,35,23,27,24,26,40 $   seems almost impossible,  unless the Susa scribes  knew this number to be a perfect square. This    might have been  revealed to them   were they capable of factoring this  number into prime numbers and observing that all the powers are even numbers\footnote{Recall that if $n=p_1^{2m_1} \times p_2^{2m_2}\times \cdots \times p_{k}^{2m_k}$, then $\sqrt{n}=p_1^{m_1} \times p_2^{m_2}\times \cdots \times p_{k}^{m_k}$.}.

\end{document}